\documentclass[12pt,leqno]{amsart}
\usepackage{amsmath,amssymb,amscd,latexsym,amsthm,mathrsfs}
\usepackage[unicode]{hyperref}
\usepackage{hypbmsec}

\ifx\pdfoutput\undefined\else
\hypersetup{
colorlinks=true,
linkcolor=blue,
citecolor=blue,
urlcolor=blue,
filecolor=blue,
bookmarksnumbered=true,
pdfstartview=FitH,
pdfhighlight=/N
}
\fi

\newtheorem{proposition}{Proposition}
\newtheorem{conjecture}[proposition]{Conjecture}
\theoremstyle{remark}

\newtheorem*{acknowledgements}{Acknowledgements}

\let\lf\lfloor
\let\rf\rfloor

\def\ord{\operatorname{ord}}

\newcommand{\qbin}[2]{\genfrac{[}{]}{0pt}{}{#1}{#2}}
\newcommand{\R}{\mathbb R}
\newcommand{\Z}{\mathbb Z}
\newcommand{\Q}{\mathbb Q}
\newcommand{\ab}{\boldsymbol{a}}
\newcommand{\bb}{\boldsymbol{b}}

\begin{document}

\title{A $q$-rious positivity}

\author{S.~Ole Warnaar}
\address{School of Mathematics and Physics, The University of Queensland,
Brisbane, QLD 4072, Australia}

\author{Wadim Zudilin}
\address{School of Mathematical and Physical Sciences,
The University of Newcastle, Callaghan, NSW 2308, Australia}

\date{}

\thanks{Work supported by the Australian Research Council.}

\begin{abstract}
The $q$-binomial coefficients $\qbin{n}{m}=
\prod_{i=1}^m (1-q^{n-m+i})/(1-q^i)$, for integers $0\leq m\leq n$,
are known to be polynomials with non-negative integer coefficients.
This readily follows from the $q$-binomial theorem, or the many 
combinatorial interpretations of $\qbin{n}{m}$.
In this note we conjecture an arithmetically motivated generalisation of
the non-negativity property for products of ratios of $q$-factorials
that happen to be polynomials.
\end{abstract}

\maketitle


The fact that the binomial coefficients
\begin{equation}
\binom{n}{m}=\frac{n!}{(n-m)!\,m!}
\label{e-bin}
\end{equation}
are integers easily follows from the following \emph{arithmetic} argument.
The order in which a prime $p$ enters~$n!$ is given by
\begin{equation}
\ord_p n!=\biggl\lf\frac{n}{p}\biggr\rf+\biggl\lf\frac{n}{p^2}\biggr\rf
+\biggl\lf\frac{n}{p^3}\biggr\rf+\dotsb,
\label{e-ord}
\end{equation}
where $\lf\,\cdot\,\rf$ is the integer-part function.
Setting $x=(n-m)/p^k$ and $y=m/p^k$ in the inequality
\[
\lf x+y\rf-\lf x\rf-\lf y\rf\ge 0,
\]
and summing $k$ over the positive integers, we see that
\[
\ord_p\binom{n}{m}\ge 0
\quad\text{for any prime $p$}.
\]
This obviously implies that $\binom{n}{m}\in\Z$.

A standard way to establish integrality purely \emph{combinatorially}
amounts to interpreting the factorial ratio in~\eqref{e-bin} as coefficients
in the expansion
\[
(1+t)^n=\sum_{m=0}^n\binom{n}{m} t^m,
\]
that is, as the number of $m$-element subsets of an $n$-set.

The arithmetic argument can be extended to more general factorial ratios.
For example, the inequality \cite[Division~8, Problems~8 and~136]{PS}
\[
\lf 2x\rf+\lf 2y\rf-\lf x\rf-\lf x+y\rf-\lf y\rf\ge0
\]
implies that
\begin{equation}
A_{n,m}:=\frac{(2n)!\,(2m)!}{n!\,(n+m)!\,m!}\in\Z.
\label{e-bin2}
\end{equation}
E.~Catalan noted this integrality as early as~1874 \cite{Catalan}.
In a much more recent study~\cite[\S\,6]{Ge},
I.~Gessel named the $A_{n,m}$ the \emph{super Catalan numbers}. He
stated several formulae for these numbers including
\begin{equation}
A_{n,m}
=\sum_{k=-\infty}^{\infty}(-1)^k\binom{2n}{n+k}\binom{2m}{m+k}.
\label{Ge1}
\end{equation}
This identity, attributed to K.~von Szily (1894)~\cite{Szily},
clearly implies the integrality claimed in \eqref{e-bin2}
but, of course, obscures that $A_{n,m}\geq 0$.

Let $\ab=(a_1,\dots,a_r)$ and $\bb=(b_1,\dots,b_s)$ be tuples of
positive integers subject to the condition
\begin{equation}
\sum_{i=1}^r\lf a_ix\rf-\sum_{j=1}^s\lf b_jx\rf\geq 0
\qquad\text{for $x\ge0$}.
\label{condb}
\end{equation}
In his work on the distribution of primes (cf.~\cite{Bo})
P.~Chebyshev considered the ratios
\begin{equation}
D_n(\ab,\bb):=\frac{(a_1n)!\dotsb(a_rn)!}{(b_1n)!\dotsb(b_sn)!}\,.
\label{e-bin3}
\end{equation}
In view of~\eqref{e-ord}, condition~\eqref{condb} is
necessary and sufficient for $D_n(\ab,\bb)\in\Z$ for all
positive~$n$\,---\,a fact known
in the literature as \emph{Landau's criterion}~\cite{La}.
Unlike the special cases of binomial coefficients and super Catalan
numbers, there seems to be no non-arithmetical approach available
in the literature to demonstrate that $D_n(\ab,\bb)\in\Z$ more generally.

It is worth mentioning that the Chebyshev--Landau factorial ratios
appear quite naturally in several deep mathematical problems including,
for example, the Riemann hypothesis~\cite{Bo}
and arithmetic properties of mirror maps~\cite{De}.
We refer the interested reader to \cite{Bo} which, among other things,
contains a full classification of pairs of tuples $(\ab,\bb)$
satisfying \eqref{condb} for $s\le r+1$ and $\sum_{i=1}^ra_i=\sum_{j=1}^sb_j$.
(The latter ``balancing'' condition in fact follows from~\eqref{condb}
if $x\geq 0$ is replaced by $x\in \R$; cf.\ \cite[Lemma~3.4]{Bo}.)

\medskip

The above integrality has an interesting $q$-counterpart.
It follows immediately from the definition of the $q$-factorial,
\[
[n]!=[n]_q!=\prod_{i=1}^n\frac{1-q^i}{1-q}\,,
\]
that $[n]!$ is a polynomial
whose irreducible factors over~$\Q$ are cyclotomic polynomials
$\Phi_\ell(q)\in\Z[q]$
(cf.~\cite[\S\,1]{Zu}). Moreover,
\[
\ord_{\Phi_\ell(q)}[n]!=\biggl\lf\frac n\ell\biggr\rf
\quad\text{for all $\ell=2,3,4,\dots$}\,.
\]
We may thus conclude that
\[
\qbin{n}{m}:=\frac{[n]!}{[n-m]!\,[m]!}\,, \qquad
A_{n,m}(q):=\frac{[2n]!\,[2m]!}{[n]!\,[n+m]!\,[m]!}\,,
\]
and
\begin{equation}
D_n(\ab,\bb;q)
:=\frac{[a_1n]!\dotsb[a_rn]!}{[b_1n]!\dotsb[b_sn]!}
\label{e-qbin3}
\end{equation}
subject to \eqref{condb} are all polynomials in $\Z[q]$.
For this reason the $q$-binomial
coefficients are often referred to by their alternative name of
\emph{Gaussian polynomials}.

Another well-known fact about the Gaussian polynomials
is the \emph{non-negativity} of their coefficients.
In fact, since each of the coefficients $c_i$ in
$\qbin{n}{m}=c_0+\cdots+c_{nm} q^{nm}$ is strictly
positive\footnote{One may in fact show that $c_0=1$ and that the
$c_i$ are symmetric and unimodal; $c_i=c_{nm-i}$ and
$c_i\leq c_{i+1}$ for $0\le i\le\lf nm/2-1\rf$.},
it is costumary to refer to them as \emph{positive polynomials}.
Following tradition, we relax the term positivity
to simply refer to any polynomial with non-negative coefficients.
Hence we say that $1+q^2$ is a positive polynomial even though the
linear term has vanishing coefficient.

The only known proofs of the positivity of the $q$-binomial coefficients
are essentially all combinatorial. For example, the $q$-binomial theorem 
\cite[Eq.~(II.4)]{GR}
\[
\prod_{i=0}^{n-1}(1+tq^i)=\sum_{m=0}^n q^{\binom{m}{2}}\qbin{n}{m}t^m,
\]
implies positivity and, more specifically, the combinatorial
interpretation
\[
q^{\binom{m}{2}}\qbin{n}{m}=
\sum_{\substack{I\subseteq \{0,\dots,n-1\} \\[1pt] |I|=m}} 
q^{\sum_{i\in I}i}.
\]

In view of the preceding discussion the following conjecture arises
naturally.
\begin{conjecture}\label{conj}
Let $\ab=(a_1,\dots,a_r)$ and $\bb=(b_1,\dots,b_s)$ satisfy \eqref{condb}.
Then the \emph{polynomial}
\begin{equation}
D(\ab,\bb;q):=\frac{[a_1]!\dotsb[a_r]!}{[b_1]!\dotsb[b_s]!}
\label{e-D}
\end{equation}
is positive.
\end{conjecture}

Replacing all $\ab$ and $\bb$ by $\ab n$ and $\bb n$ for a positive
integer~$n$, we see that the conjecture is equivalent to the claim
that the polynomials $D_n(\ab,\bb;q)$ defined in~\eqref{e-qbin3}
are positive for all positive integers~$n$.

The conjecture is trivially true whenever the right-hand side in~\eqref{e-D}
can be represented as a product of $q$-binomial coefficients. To provide
some further evidence, we show the validity of the conjecture for
the $q$-super Catalan numbers
\[
A_{n,m}(q)=D((2n,2m),(n,n+m,m);q)
\]
as well as for
\begin{equation}
B_{n,m}(q):=D((2n,m),(n,2m,n-m);q)
=\frac{[2n]![m]!}{[n]!\,[2m]!\,[n-m]!}\in\Z[q],
\quad n\ge m.
\label{e-qbin4}
\end{equation}
We note that, as shown in \cite{Bo}, the $q$-binomial coefficients
$\qbin{n}{m}$ together with $A_{n,m}(q)$ and $B_{n,m}(q)$ exhaust the
space of $2$-parameter solutions to~\eqref{condb} with $s=r+1$.

\begin{proposition}
The $q$-super Catalan numbers $A_{n,m}(q)$ are positive polynomials
for all $n,m\geq 0$.
\end{proposition}

\begin{proof}
Our proof rests on a $q$-analogue of Gessel's formula~\cite[Eq.~(32)]{Ge}
\[
A_{n,n+p}
= \sum_{k=0}^{\lfloor p/2\rfloor}
2^{p-2k} \binom{p}{2k} A_{n,k} \qquad (p\geq 0),
\]
given in \eqref{qGe3} below.

Let $n$ and $p$ be non-negative integers. Twice applying the
$q$-Chu--Vandermonde sum \cite[Eq.~(II.7)]{GR} in the form
\[
\qbin{a+b}{c}=\sum_{k=0}^{\infty} q^{k(b-c+k)} \qbin{a}{k} \qbin{b}{c-k}
\]
yields
\begin{align}
\label{e-main}
\qbin{2n+2p}{p}
&=\sum_{j=0}^{\infty} q^{j(n+j)} \qbin{n+p}{j} \qbin{n+p}{p-j}
\\
&=\sum_{j=0}^{\infty} q^{j(n+j)} \qbin{n+p}{j}
\sum_{k=0}^{\infty} q^{k(n+k)} \qbin{j}{k} \qbin{n+p-j}{p-j-k}.
\notag
\end{align}
Multiplying this by $[2n]!\,[p]!/([n]!\,[n+p]!)$ implies the recurrence
\begin{equation}
A_{n,n+p}(q)
= \sum_{k=0}^{\lfloor p/2\rfloor} A_{n,k}(q)
\sum_{j=k}^{p-k} q^{k(n+k)+j(n+j)} \qbin{p}{2k} \qbin{p-2k}{j-k}
\label{qGe3}
\end{equation}
for $p\geq 0$.
Together with the initial conditions $A_{n,n}(q)=A_{n,0}(q)=\qbin{2n}{n}$,
the symmetry $A_{n,m}(q)=A_{m,n}(q)$ and the positivity of
$q$-binomial coefficients,
formula~\eqref{qGe3} implies the desired positivity of $A_{n,m}(q)$.
\end{proof}

Another positivity result related to $A_{n,m}(q)$ may be found
in~\cite{GJZ}. Before stating this result we remark that by
taking $(a,b,c)\mapsto (1,\infty,q^{-m})$
in the very-well poised $_6\phi_5$ summation \cite[Eq.~(II.21)]{GR}
a $q$-analogue of von Szily's identity \eqref{Ge1} arises. Namely,
\[
A_{n,m}(q)=\sum_{k=-\infty}^{\infty}
(-1)^k q^{\binom{k}{2}+k^2} \qbin{2n}{n+k} \qbin{2m}{m+k}.
\]
After the substitution $q\mapsto 1/q$ this may also be written as
\begin{equation}\label{Anmq}
A_{n,m}(q)=q^{-nm}\sum_{k=-\infty}^{\infty}
(-1)^k q^{\binom{k}{2}} \qbin{2n}{n+k} \qbin{2m}{m+k}.
\end{equation}
Now a special case of \cite[Theorem 4.7]{GJZ} amounts to the following claim.
For $r,s$ positive integers and $n,m$ nonnegative integers,
the functions $R_{n,m;r,s}(q)$ defined by
\[
\sum_{k=-\infty}^{\infty}
(-1)^k q^{\binom{k}{2}} \qbin{2n}{n+k}^r \qbin{2m}{m+k}^s=
A_{n,m}(q) R_{n,m;r,s}(q)
\]
are positive polynomials.
Of course, from \eqref{Anmq} it follows that
$R_{n,m;1,1}(q)=q^{nm}$. In view of the positivity of $A_{m,n}(q)$ the
positivity of the right-hand side should come as no surprise, since,
intuitively, raising the value of $r$ and/or $s$ should
result in a ``more positive (or less-negative) polynomial''.
The fact that the right-side factors, with $A_{m,n}(q)$ as one of its
factors, is much more remarkable.
We also note that in general it is hard to prove the positivity of
alternating sum expressions of the form given above. For example,
showing that
\[
\sum_{k=-\infty}^{\infty}
(-1)^k q^{\binom{k}{2}+4k^2} \qbin{2n}{n+3k}
\]
is positive is key to proving the longstanding
Borwein conjecture, see \cite{Andrews95,Bressoud,BW,W1,W2}.

\begin{proposition}
The $B_{n,m}(q)$ defined in \eqref{e-qbin4} are positive polynomials
for all $n\geq m\geq 0$.
\end{proposition}

\begin{proof}
This time we simply multiply~\eqref{e-main} by
$[n]!\,[2n+p]!/([2n]!\,[n+p]!)$ to get
\begin{equation}
B_{n+p,n}(q)=\sum_{k=0}^{\lfloor p/2\rfloor} B_{n+k,n}(q)
\sum_{j=k}^{p-k} q^{k(n+k)+j(n+j)} \qbin{2n+p}{2n+2k} \qbin{p-2k}{j-k}.
\label{qGe4}
\end{equation}
Since $B_{n,n}(q)=1$, the result follows from~\eqref{qGe4} by
induction on~$p$.
\end{proof}

To provide additional support for Conjecture~\ref{conj} we have computed
the polynomials $D_n(\ab,\bb;q)$ for all $n$ up to $20$ for the $52$
choices for $\ab$ and $\bb$ listed in \cite[Table~2]{Bo}.
Since $s>r$ for each of these, and since the large~$n$ limit of
$D_n(\ab,\bb)$ is given by the positive power series
$\prod_{i\geq 1}(1-q^i)^{r-s}$,
one would expect potential counter examples to occur for ``small'' values
of~$n$. However, our computation resulted in polynomials with non-negative
coefficients only.

\begin{acknowledgements}
We would like to thank James Wan for pointing out the reference~\cite{Ge}.
\end{acknowledgements}

\end{document}